\documentclass[11pt]{amsart}
\usepackage{amsmath,amssymb,amsthm}
\usepackage[margin=1.15in]{geometry}
\usepackage{enumitem}
\usepackage{xcolor}
\usepackage[colorlinks=true,linkcolor=blue!50!black,citecolor=green!40!black,urlcolor=blue!50!black]{hyperref}

\newtheorem{theorem}{Theorem}[section]
\newtheorem{lemma}[theorem]{Lemma}
\newtheorem{proposition}[theorem]{Proposition}
\newtheorem{corollary}[theorem]{Corollary}
\newtheorem*{theoremHN}{Theorem HN}
\newtheorem{mainthm}{Theorem}

\theoremstyle{definition}

\newtheorem{example}[theorem]{Example}
\theoremstyle{remark}
\newtheorem{remark}[theorem]{Remark}
\newtheorem{question}[theorem]{Question}

\numberwithin{equation}{section}

\newcommand{\R}{\mathbb{R}}
\newcommand{\cB}{\mathcal{B}}
\newcommand{\cU}{\mathcal{U}}
\newcommand{\cL}{\mathcal{L}}
\newcommand{\Om}{\Omega}
\newcommand{\bOm}{\overline{\Omega}}
\newcommand{\dOm}{\partial\Omega}
\newcommand{\eps}{\varepsilon}
\newcommand{\diam}{\operatorname{diam}}
\newcommand{\dist}{\operatorname{dist}}
\newcommand{\conv}{\operatorname{conv}}
\newcommand{\norm}[1]{\|#1\|_\infty}
\newcommand{\Hf}[1]{H_{#1}}
\newcommand{\Hu}[1]{\overline{H}_{#1}}
\newcommand{\Hl}[1]{\underline{H}_{#1}}

\title[Iterated restricted means]{Convergence of iterated restricted mean values to harmonic functions}

\makeatletter
\@namedef{subjclassname@2020}{\textup{2020} Mathematics Subject Classification}
\makeatother

\author{Mohammad Javaheri}
\address{Department of Mathematics and Statistics, Siena University\\
515 Loudon Road, Loudonville, NY 12211, USA}
\email{mjavaheri@siena.edu}

\date{}


\begin{document}

\begin{abstract}
Let $\Omega\subset\R^n$ be a bounded domain and let $\delta$ be a continuous positive function on $\Omega$ with $B(x,\delta(x))\subseteq\Omega$ for every $x\in\Omega$. Let $\sigma(f)(x)$ denote the average of $f$ over the ball $B(x,\delta(x))$. Restricted mean value theorems give conditions under which a fixed point of $\sigma$ must be harmonic. We study the dynamical counterpart of this question: the behavior of the iterates $\sigma^k(f)$. We show that for \emph{every} bounded domain and every bounded measurable $f$ that is continuous at the points of $\partial\Omega$, the iterates $\sigma^k(f)$ converge locally uniformly in $\Omega$ to the Perron--Wiener--Brelot solution of the Dirichlet problem with boundary values $f|_{\partial\Omega}$, and that the convergence is uniform on $\overline\Omega$ for all continuous $f$ if and only if $\Omega$ is regular for the Dirichlet problem. This contains and extends earlier results for $C^{2,\alpha}$ domains with $C^{2,\alpha}$ data and for strictly convex domains with continuous data. We also give explicit rates of convergence: an algebraic rate $O(1/k)$ on domains with a uniform exterior ball condition and a geometric rate on convex domains, when $\delta$ is comparable to the distance to the boundary.
\end{abstract}

\maketitle

\section{Introduction}

\subsection{Restricted mean value properties}
A harmonic function on a domain $\Om\subseteq\R^n$ is equal, at each point, to its average over every ball centred at that point and contained in $\Om$. Conversely, a locally bounded measurable function with this mean value property for \emph{all} such balls is harmonic; this goes back to Koebe and is classical (see the survey of Netuka and Vesel\'y \cite{NV94}). It is natural to ask how much of the mean value property is needed. In the \emph{restricted} version, one prescribes for each point $x\in\Om$ a single radius $\delta(x)>0$ and asks only that
\begin{equation}\label{eq:rmvp}
f(x)=\frac{1}{|B_x|}\int_{B_x} f(z)\,dz,\qquad B_x:=B(x,\delta(x))\subseteq\Om,\qquad x\in\Om ,
\end{equation}
where $|E|$ denotes the Lebesgue measure of $E$ and $B(z,r)$ is the open ball of radius $r$ centred at $z$. Following \cite{HN93}, such an $f$ is called a $(\lambda,\delta)$-median function.

Volterra \cite{Vol09} and Kellogg \cite{Kel34} observed that if $f$ is continuous on $\bOm$ and satisfies \eqref{eq:rmvp}, then $f$ is harmonic; the argument is a maximum principle comparison with the solution of the Dirichlet problem. Without control up to the boundary the question becomes delicate: Littlewood's one-circle problem and its relatives were studied by Veech \cite{Vee72,Vee73} and Hansen and Nadirashvili \cite{HN93,HN94}, who found both positive results and counterexamples. We quote the following result in the form surveyed in \cite{Han96}.

\begin{theoremHN}[Hansen--Nadirashvili \cite{HN93,Han96}]
Let $\Om\subset\R^n$ be a bounded domain, $\delta$ an admissible function on $\Om$, and $f$ a $(\lambda,\delta)$-median function on $\Om$ such that $|f|\le h$ for some harmonic function $h$ on $\Om$. If $f$ is continuous, or if $\delta$ is locally bounded away from zero, then $f$ is harmonic.
\end{theoremHN}

Here and below, a function $\delta:\Om\to\R$ is called \emph{admissible} if $\delta(x)>0$ and $B(x,\delta(x))\subseteq\Om$ for all $x\in\Om$.

\subsection{Iterated averages}
Fix an admissible function $\delta$ and define the averaging operator
\begin{equation}\label{eq:sigma}
\sigma(f)(x)=
\begin{cases}
\dfrac{1}{|B_x|}\displaystyle\int_{B_x} f(z)\,dz, & x\in\Om,\\[2ex]
f(x), & x\in\dOm .
\end{cases}
\end{equation}
Restricted mean value theorems describe the \emph{fixed points} of $\sigma$. In this paper we study the \emph{dynamics} of $\sigma$: given an initial function $f$, what happens to the iterates
\[
f_0=f,\qquad f_{k+1}=\sigma(f_k),\qquad\text{i.e. } f_k=\sigma^k(f)?
\]
Since $\sigma$ fixes harmonic functions, one expects $\sigma^k(f)$ to converge to the harmonic function with the same boundary values as $f$. This is the sense in which we speak of the \emph{stability} of harmonic functions under $\sigma$.

The question has two natural sources of motivation.

\emph{Relaxation methods.} The iteration $f\mapsto\sigma(f)$ is a continuous analogue of the Jacobi (or Gauss--Seidel) relaxation scheme for the discrete Laplace equation, in which the value at each grid point is repeatedly replaced by the average of its neighbours. Convergence of $\sigma^k(f)$ shows that this ``smoothing'' recovers the solution of the Dirichlet problem, with no discretisation, for any choice of the local radius $\delta$.

\emph{Random walks.} Let $(X_k)_{k\ge0}$ be the Markov chain which, from a point $x\in\Om$, jumps to a point chosen uniformly in $B_x$ (and stays put on $\dOm$). Then $\sigma^k(f)(x)=\mathbb{E}_x[f(X_k)]$. When $\delta$ is the distance to the boundary and balls are replaced by spheres this is Muller's \emph{walk on spheres} \cite{Mul56}, a standard Monte Carlo method for the Dirichlet problem; chains of this type are also the main tool in Veech's work \cite{Vee73}. Convergence of $\sigma^k(f)$ to the solution of the Dirichlet problem says that $X_k$ converges to a boundary point distributed according to harmonic measure. The proofs below are analytic, but this interpretation explains several of the steps (Remark~\ref{rem:prob}).

\subsection{Main results}
Throughout, $\Om\subset\R^n$ ($n\ge1$) is a bounded domain, $D=\diam\Om$, $\rho(x)=\dist(x,\dOm)$, and $\delta$ is a \emph{continuous} admissible function. Since $0<\delta\le\rho$, setting $\delta=0$ on $\dOm$ gives a continuous function on $\bOm$. We write $\cB(\bOm)$ for the class of bounded functions $f:\bOm\to\R$ that are Lebesgue measurable on $\Om$ and continuous (relative to $\bOm$) at every point of $\dOm$; it contains $C(\bOm)$. For $\varphi\in C(\dOm)$, $\Hf\varphi$ denotes the Perron--Wiener--Brelot (PWB) solution of the Dirichlet problem with boundary data $\varphi$ (see \S\ref{sec:pwb}). Recall that $\Om$ is \emph{regular} if $\Hf\varphi$ extends continuously to $\bOm$ with value $\varphi$ on $\dOm$, for every $\varphi\in C(\dOm)$.

\begin{mainthm}[Arbitrary bounded domains]\label{thm:A}
Let $f\in\cB(\bOm)$ and $\varphi=f|_{\dOm}$. Then $\sigma^k(f)\to\Hf\varphi$ uniformly on compact subsets of $\Om$.
\end{mainthm}

\begin{mainthm}[Uniform convergence and regularity]\label{thm:B}
If $\Om$ is regular, then $\sigma^k(f)\to\Hf{\varphi}$ uniformly on $\bOm$ for every $f\in\cB(\bOm)$, where $\Hf\varphi$ is extended to $\bOm$ by $\varphi$. Conversely, if $\sigma^k(f)$ converges uniformly on $\bOm$ for every $f\in C(\bOm)$, then $\Om$ is regular.
\end{mainthm}

Since bounded domains satisfying an exterior cone condition, and in particular $C^{2,\alpha}$ domains and convex domains, are regular, Theorem~\ref{thm:B} contains the following corollary, which improves the two main theorems of an earlier version of this paper \cite{Jav07}: there, the first part required $f\in C^{2,\alpha}(\bOm)$ and the second required $\Om$ to be strictly convex.

\begin{corollary}\label{cor:smooth-convex}
Let $\Om$ be a bounded domain that is either of class $C^{1}$, or Lipschitz, or convex. Then for every $f\in C(\bOm)$ the iterates $\sigma^k(f)$ converge uniformly on $\bOm$ to the harmonic function $u\in C(\bOm)$ with $u=f$ on $\dOm$. If moreover $\Om$ is of class $C^{2,\alpha}$ and $f\in C^{2,\alpha}(\bOm)$ for some $0<\alpha<1$, then $u\in C^{2,\alpha}(\bOm)$.
\end{corollary}

Theorem~\ref{thm:A} also yields, as a by-product, a short proof that a bounded $(\lambda,\delta)$-median function with continuous boundary values is harmonic on an arbitrary bounded domain (Corollary~\ref{cor:rmvp}).

Finally, we quantify the convergence. Uniform convergence cannot hold at a rate that is uniform over the unit ball of $C(\bOm)$ (Remark~\ref{rem:norate}), so rates require regularity of $f$. Let $w$ denote the \emph{torsion function} of $\Om$, i.e.\ the solution of $\Delta w=-1$ in $\Om$, $w=0$ on $\dOm$.

\begin{mainthm}[Rates of convergence]\label{thm:C}
Suppose $\delta\ge\kappa\rho$ on $\Om$ for some $\kappa\in(0,1]$. Let $f\in C(\bOm)\cap C^2(\Om)$ with $\Delta f$ bounded, and let $u$ be the harmonic function on $\Om$ with $u=f$ on $\dOm$.
\begin{enumerate}[label=\textup{(\roman*)}]
\item If $\Om$ satisfies a uniform exterior ball condition with radius $R$, then
\[
\norm{\sigma^k(f)-u}\le \frac{2(n+2)A^2}{\kappa^2}\,\frac{\norm{\Delta f}}{k},\qquad A=\frac{(R+D)^{n+2}}{2R^{n+1}} .
\]
\item If $\Om$ is convex, then
\[
\norm{\sigma^k(f)-u}\le \frac{D^2}{4}\,\norm{\Delta f}\,\lambda^k,\qquad \lambda=1-\frac{\kappa^2}{16\sqrt2\,(n+2)}<1 .
\]
\end{enumerate}
\end{mainthm}

The geometric rate in (ii) is a special case of a more general estimate for convex domains (Theorem~\ref{thm:geometric}), which applies whenever $|f-u|\le M\rho^\gamma$ for some $\gamma\in(0,1)$.

\subsection{Organization}
Section~\ref{sec:prelim} collects the elementary properties of $\sigma$, including a modulus of continuity (Lemma~\ref{lem:modulus}) and an ``energy identity'' (Lemma~\ref{lem:energy}) which is the key to Theorem~\ref{thm:A}, and recalls the facts we need from potential theory. Theorem~\ref{thm:A} is proved in Section~\ref{sec:general}, Theorem~\ref{thm:B} and its corollaries in Section~\ref{sec:regular}, and Theorem~\ref{thm:C} in Section~\ref{sec:rates}. Section~\ref{sec:convexhull} presents an alternative, geometric argument based on convex hulls of graphs, which was the approach of \cite{Jav07}. Section~\ref{sec:remarks} contains further remarks and open questions.

\subsection*{Notation}
$\omega_n$ is the volume of the unit ball of $\R^n$; $\delta_x=\delta(x)$; $\norm{f}=\sup_{\bOm}|f|$ for bounded $f$ (all functions considered are bounded, and $\sigma(f)$ does not change if $f$ is modified on a null subset of $\Om$, so nothing is lost by using the supremum instead of the essential supremum). For $y\in\R^n$ we write $q_y(z)=|z-y|^2$.

\section{Preliminaries}\label{sec:prelim}

\subsection{Elementary properties of $\sigma$}

\begin{lemma}\label{lem:basic}
For bounded functions $f,g$ on $\bOm$ that are measurable on $\Om$:
\begin{enumerate}[label=\textup{(\alph*)}]
\item $\sigma$ is linear, positive ($f\le g\Rightarrow\sigma f\le\sigma g$) and $\sigma(1)=1$; hence $\norm{\sigma f}\le\norm f$ and $\norm{\sigma f-\sigma g}\le\norm{f-g}$.
\item For $x\in\Om$, $\sigma^k(f)(x)$ depends only on the restriction of $f$ to $\Om$.
\item For each $x\in\bOm$ and $k\ge0$ there is a Borel probability measure $\nu^{(k)}_x$ on $\bOm$ with $\sigma^k(f)(x)=\int f\,d\nu_x^{(k)}$. In particular, Jensen's inequality $\sigma^k(f^2)\ge(\sigma^kf)^2$ holds.
\end{enumerate}
\end{lemma}

\begin{proof}
(a) and (b) are immediate from \eqref{eq:sigma}; for (b) note that $B_x\subseteq\Om$ for $x\in\Om$. For (c), let $\nu_x$ be normalized Lebesgue measure on $B_x$ when $x\in\Om$ and the Dirac mass at $x$ when $x\in\dOm$. By Lemma~\ref{lem:modulus} below, $x\mapsto\nu_x(E)$ is continuous on $\Om$ for every Borel set $E$, so $(\nu_x)$ is a Markov kernel and $\nu^{(k)}_x$ is its $k$-fold composition.
\end{proof}

\begin{lemma}[Modulus of continuity]\label{lem:modulus}
Let $f$ be bounded on $\Om$ and measurable. Then for all $x,y\in\Om$,
\begin{equation}\label{eq:modulus}
|\sigma(f)(x)-\sigma(f)(y)|\le n2^{n+1}\,\norm f\,\frac{|\delta_x-\delta_y|+|x-y|}{\delta_x}.
\end{equation}
In particular $\sigma(f)$ is continuous on $\Om$, and if $\delta$ is Lipschitz then $\sigma(f)$ is locally Lipschitz.
\end{lemma}

\begin{proof}
Put $\eta=|\delta_x-\delta_y|+|x-y|$. If $\eta>\delta_x$, then \eqref{eq:modulus} follows from $|\sigma f(x)-\sigma f(y)|\le2\norm f$. Suppose $\eta\le\delta_x$. If $|z-x|<\delta_x-\eta$ then $|z-y|<\delta_x-\eta+|x-y|\le\delta_y$, and if $|z-y|<\delta_y$ then $|z-x|<\delta_y+|x-y|\le\delta_x+\eta$. Hence both $B_x$ and $B_y$ contain $B(x,\delta_x-\eta)$ and are contained in $B(x,\delta_x+\eta)$, so that
\[
|B_x\triangle B_y|\le\omega_n\big((\delta_x+\eta)^n-(\delta_x-\eta)^n\big)\le 2n\omega_n\eta(\delta_x+\eta)^{n-1}\le n2^n\omega_n\eta\,\delta_x^{n-1}.
\]
Now
\[
\sigma f(x)-\sigma f(y)=\frac{1}{|B_x|}\Big(\int_{B_x}f-\int_{B_y}f\Big)+\Big(\frac{1}{|B_x|}-\frac1{|B_y|}\Big)\int_{B_y}f,
\]
and each term is bounded in absolute value by $\norm f\,|B_x\triangle B_y|/|B_x|$, since $\big||B_x|-|B_y|\big|\le|B_x\triangle B_y|$. Combining, $|\sigma f(x)-\sigma f(y)|\le 2\norm f\, n2^n\eta/\delta_x$.
\end{proof}

\begin{lemma}[Boundary behaviour]\label{lem:boundary}
Let $f$ be bounded on $\bOm$, measurable on $\Om$, let $y\in\dOm$, and suppose $f(x)\to L$ as $x\to y$, $x\in\Om$. Then $\sigma(f)(x)\to L$ as $x\to y$, $x\in\Om$. Consequently $\sigma$ maps $\cB(\bOm)$ into $\cB(\bOm)\cap C(\Om)$ and $C(\bOm)$ into $C(\bOm)$.
\end{lemma}

\begin{proof}
Given $a>0$ choose $b>0$ with $|f(x)-L|\le a$ whenever $x\in\Om$ and $|x-y|<b$. If $z\in\Om$ and $|z-y|<b/2$, then $\delta_z\le\rho(z)\le|z-y|<b/2$, so every $x\in B_z$ satisfies $|x-y|<b$, and therefore $|\sigma f(z)-L|\le|B_z|^{-1}\int_{B_z}|f-L|\le a$. The last statement follows from this, Lemma~\ref{lem:modulus}, and the fact that $\sigma f=f$ on $\dOm$.
\end{proof}

\subsection{Test functions}
The following identities are the only place where the specific geometry of balls enters.

\begin{lemma}\label{lem:test}
\begin{enumerate}[label=\textup{(\alph*)}]
\item If $u$ is bounded and harmonic on $\Om$, then $\sigma u=u$ on $\Om$.
\item If $s$ is superharmonic and bounded below on $\Om$, then $\sigma s\le s$ on $\Om$.
\item For every $y\in\R^n$, $\sigma(q_y)=q_y+\frac{n}{n+2}\delta^2$ on $\bOm$.
\item If $w\in C^2(\Om)$ is bounded and $\Delta w=-1$ in $\Om$, then $\sigma w=w-\frac{1}{2(n+2)}\delta^2$ on $\Om$.
\end{enumerate}
\end{lemma}

\begin{proof}
(a), (b): For $r<\delta_x$ the closed ball $\overline B(x,r)$ lies in $\Om$, so the mean value of $u$ (resp.\ $s$) over $B(x,r)$ equals (resp.\ is at most) its value at $x$. Letting $r\uparrow\delta_x$ and using dominated (resp.\ monotone) convergence gives the claim; note that $B_x$ may touch $\dOm$, which is why boundedness is assumed.
(c): For $x\in\Om$, expanding $|z-y|^2=|x-y|^2+2(z-x)\cdot(x-y)+|z-x|^2$ and using $\int_{B(x,r)}(z-x)\,dz=0$ and $\int_{B(0,r)}|z|^2dz=\frac{n}{n+2}\omega_nr^{n+2}$ gives $\sigma(q_y)(x)=q_y(x)+\frac{n}{n+2}\delta_x^2$. On $\dOm$ both sides equal $q_y$.
(d): $w+q_0/(2n)$ is bounded and harmonic, so by (a) and (c), $\sigma w=w+\frac{1}{2n}(q_0-\sigma q_0)=w-\frac{1}{2(n+2)}\delta^2$.
\end{proof}

\subsection{An energy identity}
The next lemma says, in probabilistic language, that the chain $(X_k)$ visits any compact subset of $\Om$ only finitely often on average.

\begin{lemma}[Energy identity]\label{lem:energy}
For every $x\in\bOm$,
\[
\sum_{k=0}^\infty\sigma^k(\delta^2)(x)\le\frac{n+2}{n}\,\sup_{z\in\bOm}|z-x|^2\le\frac{n+2}{n}D^2 .
\]
Consequently, if $K\subset\Om$ is compact and $\beta_K=\min_K\delta>0$, then
\begin{equation}\label{eq:visits}
\sum_{k=0}^{\infty}\sigma^k(\mathbf 1_K)\le\frac{(n+2)D^2}{n\,\beta_K^2}\quad\text{on }\bOm;\qquad\text{in particular } \sigma^k(\mathbf 1_K)\to0 \text{ pointwise on }\bOm.
\end{equation}
\end{lemma}

\begin{proof}
By Lemma~\ref{lem:test}(c) and linearity, $\sigma^{k+1}q_x-\sigma^kq_x=\sigma^k(\sigma q_x-q_x)=\frac{n}{n+2}\sigma^k(\delta^2)$. Summing over $0\le k<N$ and evaluating at $x$, where $q_x(x)=0$,
\[
\frac{n}{n+2}\sum_{k=0}^{N-1}\sigma^k(\delta^2)(x)=\sigma^Nq_x(x)\le\sup_{\bOm}q_x .
\]
The second statement follows from $\mathbf 1_K\le\delta^2/\beta_K^2$ and the positivity of $\sigma$.
\end{proof}

\begin{remark}[Probabilistic interpretation]\label{rem:prob}
With $(X_k)$ the chain described in the introduction, each coordinate of $X_k$ is a bounded martingale, so $X_k$ converges almost surely to some $X_\infty$. Lemma~\ref{lem:energy} is the statement $\mathbb E_x\sum_k|X_{k+1}-X_k|^2\le D^2$, which forces $\delta(X_k)\to0$ and hence $X_\infty\in\dOm$; the bound \eqref{eq:visits} controls the expected number of visits to $K$. Theorem~\ref{thm:A} then says that $X_\infty$ has the harmonic measure of $\Om$ relative to $x$ as its law. In particular, the chain never ``sees'' polar parts of the boundary, such as the centre of a punctured disc.
\end{remark}

\subsection{Potential theory}\label{sec:pwb}
We recall the facts we use; see \cite[Ch.~6]{AG01}, \cite{Doo84} or \cite[\S2.8]{GT83}. Let $\varphi\in C(\dOm)$. The \emph{upper class} $\cU_\varphi$ consists of the superharmonic functions $s$ on $\Om$ that are bounded below and satisfy $\liminf_{x\to y}s(x)\ge\varphi(y)$ for every $y\in\dOm$; the lower class is $\cL_\varphi=-\cU_{-\varphi}$. Put $\Hu\varphi=\inf\cU_\varphi$ and $\Hl\varphi=\sup\cL_\varphi$. By Wiener's resolutivity theorem, $\Hu\varphi=\Hl\varphi=:\Hf\varphi$ and $\Hf\varphi$ is a bounded harmonic function on $\Om$, the PWB solution; moreover $\min\varphi\le \Hf\varphi\le\max\varphi$.

A point $y\in\dOm$ is \emph{regular} if $\Hf\varphi(x)\to\varphi(y)$ as $x\to y$ for every $\varphi\in C(\dOm)$, and $\Om$ is regular if all its boundary points are. All boundary points of bounded intervals are regular; so are those satisfying an exterior cone condition (Poincar\'e--Zaremba), in particular all boundary points of Lipschitz, $C^1$ and convex domains; in the plane every bounded simply connected domain is regular. The centre of a punctured ball and the tip of Lebesgue's spike are irregular. Regularity is characterized by Wiener's criterion.

If $\Om$ is regular, the \emph{torsion function} $w=\Hf\psi-\psi$, where $\psi(x)=|x|^2/(2n)$, belongs to $C(\bOm)\cap C^\infty(\Om)$ and satisfies $\Delta w=-1$ in $\Om$ and $w=0$ on $\dOm$. By the minimum principle for the superharmonic function $w$, $w>0$ in $\Om$.

\section{Convergence in arbitrary bounded domains}\label{sec:general}

\begin{proof}[Proof of Theorem~\ref{thm:A}]
Let $f\in\cB(\bOm)$ and $\varphi=f|_{\dOm}\in C(\dOm)$.

\emph{Step 1: $\limsup_k\sigma^kf\le\Hu\varphi$ on $\Om$.} Let $s\in\cU_\varphi$. Replacing $s$ by $\min\{s,\norm f+1\}$, which is still in $\cU_\varphi$, we may assume that $s$ is bounded. Fix $\eps>0$. For each $y\in\dOm$ there is $r_y>0$ such that $s(x)>\varphi(y)-\eps/2$ and $f(x)<\varphi(y)+\eps/2$ for all $x\in\Om\cap B(y,r_y)$; here we used the definition of $\cU_\varphi$ and the continuity of $f$ at $y$. Cover $\dOm$ by finitely many of these balls and let $N$ be their union. Then $K=\bOm\setminus N$ is a compact subset of $\Om$, and $f-s<\eps$ on $\Om\setminus K$. With $M=\sup_K(f-s)^+<\infty$ we obtain
\[
f\le s+\eps+M\mathbf 1_K\qquad\text{on }\Om .
\]
Fix $x\in\Om$. By Lemma~\ref{lem:basic}(b) only values on $\Om$ matter; by positivity and Lemma~\ref{lem:test}(b), $\sigma^k s\le\sigma^{k-1}s\le\dots\le s$ on $\Om$. Hence
\[
\sigma^kf(x)\le s(x)+\eps+M\sigma^k(\mathbf 1_K)(x),
\]
and by \eqref{eq:visits}, $\limsup_k\sigma^kf(x)\le s(x)+\eps$. Since $\eps>0$ and $s\in\cU_\varphi$ are arbitrary, $\limsup_k\sigma^kf\le\Hu\varphi$.

\emph{Step 2: pointwise convergence.} Applying Step~1 to $-f$ gives $\liminf_k\sigma^kf\ge\Hl\varphi$. By Wiener's theorem $\Hu\varphi=\Hl\varphi=\Hf\varphi$, so $\sigma^kf\to\Hf\varphi$ pointwise on $\Om$.

\emph{Step 3: local uniformity.} Let $K'\subset\Om$ be compact, and let $\beta=\min_{K'}\delta>0$. For $k\ge1$, Lemma~\ref{lem:modulus} applied to $\sigma^{k-1}f$, together with $\norm{\sigma^{k-1}f}\le\norm f$, gives
\[
|\sigma^kf(x)-\sigma^kf(y)|\le n2^{n+1}\norm f\,\beta^{-1}\big(|\delta_x-\delta_y|+|x-y|\big),\qquad x,y\in K'.
\]
Since $\delta$ is uniformly continuous on $K'$, the family $\{\sigma^kf\}_{k\ge1}$ is equicontinuous on $K'$. An equicontinuous sequence converging pointwise on a compact set converges uniformly there.
\end{proof}

\begin{remark}
The limit does not depend on $\delta$: every continuous admissible averaging scheme recovers the same function $\Hf\varphi$. Also, only the boundary values of $f$ and the values of $f$ on $\Om$ up to null sets enter; $f$ need not be continuous in $\Om$. After one step, $\sigma f$ is continuous in $\Om$, so the iteration is regularizing.
\end{remark}

As an immediate consequence we obtain a restricted mean value theorem valid on every bounded domain.

\begin{corollary}\label{cor:rmvp}
Let $f\in\cB(\bOm)$ satisfy \eqref{eq:rmvp}. Then $f=\Hf{f|_{\dOm}}$ on $\Om$; in particular $f$ is harmonic in $\Om$.
\end{corollary}

\begin{proof}
$\sigma^kf=f$ on $\Om$ for all $k$, and Theorem~\ref{thm:A} applies.
\end{proof}

Corollary~\ref{cor:rmvp} is weaker than Theorem~HN, which needs no information at the boundary, but its proof is elementary given the PWB method, and it does not require $\Om$ to be regular, unlike the classical maximum principle argument of Volterra and Kellogg.

\section{Uniform convergence and regular domains}\label{sec:regular}

\begin{proof}[Proof of Theorem~\ref{thm:B}]
Suppose first that $\Om$ is regular, let $f\in\cB(\bOm)$, $\varphi=f|_{\dOm}$, and let $u\in C(\bOm)$ be $\Hf\varphi$ extended by $\varphi$. By Lemma~\ref{lem:test}(a), $\sigma u=u$ on $\bOm$, hence $\sigma^kf-u=\sigma^k g$ with $g=f-u$. The function $g$ is bounded, measurable on $\Om$, vanishes on $\dOm$ and is continuous at each boundary point.

\emph{Step 1: $\norm{\sigma^kw}\to0$}, where $w$ is the torsion function. By Lemma~\ref{lem:test}(d), $\sigma w=w-\frac{1}{2(n+2)}\delta^2\le w$ on $\bOm$, so by positivity $(\sigma^kw)_k$ is a non-increasing sequence of non-negative functions, each continuous on $\bOm$ by Lemma~\ref{lem:boundary}. By Theorem~\ref{thm:A} it converges pointwise on $\Om$ to $\Hf0=0$, and it vanishes on $\dOm$. By Dini's theorem the convergence is uniform on $\bOm$.

\emph{Step 2.} Let $\eps>0$. Since $g$ is continuous with value $0$ at every point of the compact set $\dOm$, there is a compact $K\subset\Om$ with $|g|<\eps$ on $\bOm\setminus K$. Let $c_K=\min_Kw>0$. Then
\[
|g|\le\eps+\frac{\norm g}{c_K}\,w\qquad\text{on }\bOm,
\]
and by positivity, $\sigma(1)=1$ and Step~1,
\[
\limsup_{k\to\infty}\norm{\sigma^kg}\le\eps+\frac{\norm g}{c_K}\lim_{k\to\infty}\norm{\sigma^kw}=\eps .
\]
Hence $\sigma^kf\to u$ uniformly on $\bOm$.

Conversely, suppose that $\sigma^kf$ converges uniformly on $\bOm$ for every $f\in C(\bOm)$. Let $\varphi\in C(\dOm)$ and extend it to $f\in C(\bOm)$ by Tietze's theorem. The iterates are continuous on $\bOm$ (Lemma~\ref{lem:boundary}), so their uniform limit $v$ is continuous on $\bOm$. On $\dOm$, $\sigma^kf=f$, so $v=\varphi$ there; on $\Om$, $v=\Hf\varphi$ by Theorem~\ref{thm:A}. Thus $\Hf\varphi(x)\to\varphi(y)$ as $x\to y\in\dOm$, for every $\varphi$, i.e.\ $\Om$ is regular.
\end{proof}

\begin{proof}[Proof of Corollary~\ref{cor:smooth-convex}]
Lipschitz and $C^1$ domains satisfy an exterior cone condition at every boundary point. A convex domain has a supporting half-space at each boundary point, whose complement is an exterior cone (indeed an exterior ball). Hence these domains are regular and the first statement follows from Theorem~\ref{thm:B}. The second statement is Kellogg's theorem \cite[Thm.~6.14]{GT83}, which requires $0<\alpha<1$.
\end{proof}

\begin{remark}[Comparison with the earlier version]
In \cite{Jav07} the first part of Corollary~\ref{cor:smooth-convex} was stated for $C^{2,\alpha}$ domains and data with $\alpha\in[0,1]$. Schauder theory, and hence the existence of $u\in C^{2,\alpha}(\bOm)$, fails for $\alpha\in\{0,1\}$, but Theorem~\ref{thm:B} shows that no smoothness of $f$ is needed for the convergence itself. Similarly, strict convexity, used in \cite{Jav07} through the convex hull argument of Section~\ref{sec:convexhull}, is unnecessary.
\end{remark}

\begin{remark}[An alternative proof for continuous data]\label{rem:poly}
For $f\in C(\bOm)$, Step~2 can be replaced by polynomial approximation. If $p$ is a polynomial, then $p-\Hf p$ satisfies $\Delta(p-\Hf p)=\Delta p$, and comparing with $\pm\norm{\Delta p}w$ via the maximum principle gives $|p-\Hf p|\le\norm{\Delta p}\,w$ (Lemma~\ref{lem:torsion}(a)). Hence $\norm{\sigma^kp-\Hf p}\le\norm{\Delta p}\norm{\sigma^kw}\to0$. For general $f$, choose $p$ with $\norm{f-p}<\eps$ (Stone--Weierstrass) and use $\norm{\sigma^k(f-p)}\le\eps$ and $\norm{\Hf{f}-\Hf p}\le\eps$. Thus on a regular domain everything reduces to the single fact $\norm{\sigma^kw}\to0$.
\end{remark}

\begin{example}[An irregular domain]\label{ex:punctured}
Let $\Om=B(0,1)\setminus\{0\}\subset\R^2$ and $f(x)=1-|x|\in C(\bOm)$. The origin is polar, so $\Hf{f|_{\dOm}}$ is the harmonic extension of the zero function on the unit circle, namely $0$. By Theorem~\ref{thm:A}, $\sigma^kf\to0$ locally uniformly in $\Om$. On the other hand $\sigma^kf(0)=1$ for all $k$, and by Lemma~\ref{lem:boundary} each $\sigma^kf$ is continuous at $0$, so $\norm{\sigma^kf}=1$ for every $k$. The convergence is not uniform on $\bOm$, as Theorem~\ref{thm:B} predicts. In terms of Remark~\ref{rem:prob}: started near the origin, the chain needs many steps to escape (since $\log|X_k|$ is a martingale with bounded increments), but it eventually leaves and lands on the unit circle.
\end{example}

\section{Rates of convergence}\label{sec:rates}

\begin{remark}[No uniform rate]\label{rem:norate}
Let $C_0=\{g\in C(\bOm):g|_{\dOm}=0\}$. The restriction of $\sigma^k$ to $C_0$ has norm $1$ for every $k$: take $g_j\in C_0$ with $0\le g_j\le1$ and $g_j\uparrow\mathbf 1_\Om$; by monotone convergence $\sigma^kg_j\uparrow\sigma^k(\mathbf 1_\Om)=1$ on $\Om$. Thus no rate of convergence can hold uniformly over bounded sets of $C(\bOm)$, and rates must depend on the regularity of $f$ relative to the boundary.
\end{remark}

\begin{lemma}\label{lem:torsion}
Let $\Om$ be regular with torsion function $w$.
\begin{enumerate}[label=\textup{(\alph*)}]
\item If $f\in C(\bOm)\cap C^2(\Om)$ and $\Delta f$ is bounded, and $u$ is harmonic with $u=f$ on $\dOm$, then $|f-u|\le\norm{\Delta f}\,w$ on $\bOm$.
\item If $\Om$ is convex, then $w\le\frac D2\rho$.
\item If $\Om$ satisfies a uniform exterior ball condition with radius $R$, i.e.\ for every $y\in\dOm$ there is a ball $B(c,R)$ with $B(c,R)\cap\Om=\emptyset$ and $y\in\partial B(c,R)$, then $w\le A\rho$ with $A=(R+D)^{n+2}/(2R^{n+1})$.
\end{enumerate}
\end{lemma}

\begin{proof}
(a) With $M=\norm{\Delta f}$, the functions $Mw\pm(f-u)$ are continuous on $\bOm$, vanish on $\dOm$ and satisfy $\Delta(Mw\pm(f-u))=-M\pm\Delta f\le0$, so they are non-negative by the minimum principle.

(b) Fix $x\in\Om$ and a nearest point $y\in\dOm$. The ball $B(x,\rho(x))\subseteq\Om$ touches $\dOm$ at $y$, so any supporting hyperplane of the convex set $\Om$ at $y$ is tangent to this ball; thus $\Om\subseteq H=\{z:(z-y)\cdot e>0\}$ with $e=(x-y)/|x-y|$. Put $t(z)=(z-y)\cdot e$; then $0\le t\le D$ on $\bOm$, and $v=t(D-t)/2$ satisfies $\Delta v=-1$ and $v\ge0$ on $\dOm$. Since $v-w$ is harmonic and non-negative on $\dOm$, $w\le v$, and at $z=x$, where $t(x)=\rho(x)$, we get $w(x)\le D\rho(x)/2$.

(c) Fix $x\in\Om$, a nearest point $y\in\dOm$, and an exterior ball $B(c,R)$ at $y$. Let $b(z)=K\big(R^{-n}-|z-c|^{-n}\big)$ with $K=(R+D)^{n+2}/(2n)$. Then $|z-c|\ge R$ on $\bOm$, so $b\ge0$ there, and since $\Delta(|z|^{-n})=2n|z|^{-n-2}$ and $|z-c|\le D+R$ on $\bOm$,
\[
\Delta b(z)=-2nK|z-c|^{-n-2}\le-1\qquad(z\in\Om).
\]
So $b-w$ is superharmonic and non-negative on $\dOm$, whence $w\le b$. Finally $|x-c|\le|x-y|+R=\rho(x)+R$, and the mean value theorem gives $b(x)\le KnR^{-n-1}(|x-c|-R)\le KnR^{-n-1}\rho(x)=A\rho(x)$.
\end{proof}

Every bounded $C^{1,1}$ domain, and in particular every bounded $C^2$ domain, satisfies a uniform exterior ball condition.

\begin{proof}[Proof of Theorem~\ref{thm:C}\textup{(i)}]
By Lemma~\ref{lem:torsion}(a) and positivity, $|\sigma^kf-u|=|\sigma^k(f-u)|\le\norm{\Delta f}\,\sigma^kw$, so it suffices to show $\norm{\sigma^kw}\le 1/(ck)$ with $c=\kappa^2/(2(n+2)A^2)$. By Lemma~\ref{lem:test}(d), $\delta\ge\kappa\rho$ and Lemma~\ref{lem:torsion}(c),
\[
\sigma w=w-\frac{\delta^2}{2(n+2)}\le w-\frac{\kappa^2\rho^2}{2(n+2)}\le w-c\,w^2 .
\]
Applying $\sigma^k$ and Jensen's inequality (Lemma~\ref{lem:basic}(c)),
\[
a_{k+1}\le a_k-c\,a_k^2,\qquad a_k:=\sigma^kw(x),
\]
for each fixed $x\in\bOm$. If $a_{k+1}>0$ then $ca_k<1$ and $1/a_{k+1}\ge 1/(a_k(1-ca_k))\ge1/a_k+c$; by induction $a_k\le1/(ck)$ (trivially so once some $a_j=0$).
\end{proof}

For convex domains we obtain a geometric rate by using a power of the distance function as a Lyapunov function.

\begin{theorem}\label{thm:geometric}
Let $\Om$ be a bounded convex domain, $\delta\ge\kappa\rho$ with $\kappa\in(0,1]$, and $\gamma\in(0,1)$. Then
\[
\sigma(\rho^\gamma)\le\lambda_\gamma\,\rho^\gamma\quad\text{on }\bOm,\qquad \lambda_\gamma:=1-\frac{\gamma(1-\gamma)\,2^{\gamma-3}\kappa^2}{n+2}<1 .
\]
Consequently, if $f\in\cB(\bOm)$, $u$ is the harmonic function with $u=f$ on $\dOm$, and $|f-u|\le M\rho^\gamma$ on $\Om$, then
\[
|\sigma^kf-u|\le M\lambda_\gamma^k\rho^\gamma\le M\Big(\frac D2\Big)^\gamma\lambda_\gamma^k\qquad\text{on }\bOm.
\]
\end{theorem}

\begin{proof}
Fix $x\in\Om$ and let $y$, $e$, $H$ and $t$ be as in the proof of Lemma~\ref{lem:torsion}(b). Since $\Om\subseteq H$, the segment from $z\in\Om$ to its orthogonal projection on $\partial H$ leaves $\Om$, so $\rho\le t$ on $\Om$, with equality at $x$. Writing $z=x+\delta_x\xi$ with $\xi$ uniform in the unit ball, and $Z=\xi\cdot e$, $\tau=\delta_x/\rho(x)\in[\kappa,1]$,
\[
\sigma(\rho^\gamma)(x)\le\frac1{|B_x|}\int_{B_x}t^\gamma=\rho(x)^\gamma\,\Phi(\tau),\qquad \Phi(\tau):=\mathbb E\big[(1+\tau Z)^\gamma\big].
\]
Since $Z$ is symmetric and $s\mapsto(1+s)^{\gamma-1}$ is decreasing, $\Phi'(\tau)=\gamma\,\mathbb E[Z(1+\tau Z)^{\gamma-1}]\le0$, so $\Phi(\tau)\le\Phi(\kappa)$. For $s\in(-1,1]$, Taylor's theorem and $(1+\xi)^{\gamma-2}\ge2^{\gamma-2}$ for $\xi\le1$ give $(1+s)^\gamma\le1+\gamma s-\gamma(1-\gamma)2^{\gamma-3}s^2$. With $\mathbb EZ=0$ and $\mathbb EZ^2=1/(n+2)$ we get $\Phi(\kappa)\le\lambda_\gamma$. On $\dOm$ both sides vanish.

For the second statement, $\sigma u=u$ and positivity give $|\sigma^kf-u|=|\sigma^k(f-u)|\le M\sigma^k(\rho^\gamma)\le M\lambda_\gamma^k\rho^\gamma$, and $\rho\le D/2$.
\end{proof}

\begin{proof}[Proof of Theorem~\ref{thm:C}\textup{(ii)}]
By Lemma~\ref{lem:torsion}(a),(b) and $\rho\le D/2$,
\[
|f-u|\le\norm{\Delta f}\,\tfrac D2\,\rho\le\norm{\Delta f}\big(\tfrac D2\big)^{3/2}\rho^{1/2}.
\] Apply Theorem~\ref{thm:geometric} with $\gamma=\frac12$, for which $\lambda_{1/2}=1-\kappa^2/(16\sqrt2(n+2))$.
\end{proof}

\begin{remark}
The constants are not optimized. Numerically, for $\Om=(0,1)$ and $\delta=\rho$, one finds $\norm{\sigma^kw}\approx 6.5\times10^{-4}$ for $k=50$ and $\approx1.4\times10^{-5}$ for $k=100$, a contraction factor of roughly $0.93$ per step, compared with the bound $\lambda_{1/2}\approx0.9926$. The hypothesis $\delta\ge\kappa\rho$ cannot simply be dropped: if $\delta$ is much smaller than $\rho$ on a large set, the iteration moves mass very slowly.
\end{remark}

\section{A geometric argument: convex hulls of graphs}\label{sec:convexhull}

This section presents the argument of \cite{Jav07}, in corrected and slightly more general form. It is independent of potential theory and gives a direct reason for equicontinuity at the boundary. For $f\in C(\bOm)$ let
\[
\mathcal H(f)=\conv\{(x,f(x)):x\in\bOm\}\subset\R^{n+1}
\]
be the convex hull of the graph of $f$. By Carath\'eodory's theorem \cite{Eck93},
\[
\mathcal H(f)=\Big\{\sum_{i=1}^{n+2}\alpha_i\big(x_i,f(x_i)\big):\ \alpha_i\ge0,\ \sum_i\alpha_i=1,\ x_i\in\bOm\Big\},
\]
and it is compact.

\begin{proposition}\label{prop:hull}
Let $f\in C(\bOm)$. Then $\mathcal H(\sigma f)\subseteq\mathcal H(f)$, hence $(x,\sigma^kf(x))\in\mathcal H(f)$ for all $x\in\bOm$ and $k\ge0$.
\end{proposition}

\begin{proof}
It suffices to show $(x,\sigma f(x))\in\mathcal H(f)$. For $x\in\dOm$ this is clear. For $x\in\Om$, the point $(x,\sigma f(x))$ is the barycentre of the probability measure obtained by pushing normalized Lebesgue measure on $B_x$ forward under $z\mapsto(z,f(z))$, because the barycentre of $B_x$ is $x$. The barycentre of a probability measure on a compact set lies in the closed convex hull of that set (by the Hahn--Banach separation theorem, or by approximating the integral with Riemann sums as in \cite{Jav07}).
\end{proof}

\begin{proposition}\label{prop:extreme}
Let $f\in C(\bOm)$ and let $y\in\dOm$ be an extreme point of $\conv\bOm$. Then for every $a>0$ there is $b>0$ such that $|t-f(y)|<a$ whenever $(z,t)\in\mathcal H(f)$ and $|z-y|<b$. In particular, the iterates $\{\sigma^kf\}_{k\ge0}$ are equicontinuous at $y$.
\end{proposition}

\begin{proof}
If not, there exist $(z_j,t_j)\in\mathcal H(f)$ with $z_j\to y$ and $|t_j-f(y)|\ge a$. By compactness we may assume $t_j\to t$, so $(y,t)\in\mathcal H(f)$ and $|t-f(y)|\ge a$. Write $(y,t)=\sum\alpha_i(x_i,f(x_i))$ with $x_i\in\bOm$. Then $y=\sum\alpha_ix_i$, and since $y$ is extreme in $\conv\bOm$, $x_i=y$ whenever $\alpha_i>0$. Hence $t=f(y)$, a contradiction. The last assertion follows from Proposition~\ref{prop:hull}.
\end{proof}

If $\Om$ is strictly convex, i.e.\ every non-trivial convex combination of points of $\bOm$ lies in $\Om$, then every point of $\dOm$ is extreme and Proposition~\ref{prop:extreme} yields equicontinuity on all of $\bOm$ (together with Lemma~\ref{lem:modulus} in the interior). One can then conclude with the Arzel\`a--Ascoli theorem and the following elementary maximum principle, which replaces Theorem~\ref{thm:A} in this setting.

\begin{lemma}\label{lem:maxprinciple}
Let $g\in C(\bOm)$ with $g=0$ on $\dOm$ and $\sigma g=g$. Then $g\equiv0$.
\end{lemma}

\begin{proof}
Suppose $m=\max g>0$. The set $Z=\{g=m\}$ is a non-empty compact subset of $\Om$; choose $x_0\in Z$ maximizing the first coordinate. Since $g\le m$ is continuous and its mean over $B_{x_0}$ equals $m$, $g=m$ on $B_{x_0}$, which contains points of $Z$ with larger first coordinate. Hence $\max g\le0$; likewise $\min g\ge0$.
\end{proof}

Indeed, assume $f=0$ on $\dOm$ (subtract the harmonic extension), let $\alpha=\lim_k\norm{\sigma^kf}$ (the sequence is non-increasing), and let $F$ be a uniform limit of a subsequence $\sigma^{k_j}f$. Then $F\in C(\bOm)$, $F=0$ on $\dOm$, and $\norm{\sigma^mF}=\alpha$ for every $m$, because $\sigma^mF=\lim_j\sigma^{k_j+m}f$ uniformly. Arguing as in Lemma~\ref{lem:maxprinciple} with the nested compact sets $\Lambda_i^\pm=\{\pm\sigma^iF=\alpha\}$, one of which is non-empty for infinitely many, hence (by nestedness) all, $i$, shows that $\alpha=0$.

\begin{remark}
In \cite{Jav07}, the sets $\{\sigma^iF=\alpha\}$ were used without the sign $\pm$; since $\norm{\cdot}$ is a supremum of absolute values, both signs must be considered. Also, the barrier there was justified by the concavity of the function $h$ solving $\Delta h=-1$, $h|_{\dOm}=0$; such $h$ need not be concave, and the correct reason is that it is superharmonic (Lemma~\ref{lem:test}(d)), indeed $\sigma h=h-\delta^2/(2(n+2))$.
\end{remark}

\section{Further remarks and questions}\label{sec:remarks}

\begin{remark}[Spherical means]
Let $\sigma_S(f)(x)$ be the average of $f$ over the sphere $\partial B_x$, and assume $\overline{B_x}\subset\Om$ for $x\in\Om$ (i.e.\ $\delta<\rho$). For $f\in C(\bOm)$ all the ingredients used above remain valid: $\sigma_S$ maps $C(\bOm)$ into itself, harmonic functions are fixed, superharmonic functions satisfy $\sigma_Ss\le s$, and $\sigma_S(q_y)=q_y+\delta^2$, which gives the energy identity. Hence Theorem~\ref{thm:A} holds pointwise in $\Om$ for $f\in C(\bOm)$, and Theorem~\ref{thm:B} holds verbatim for $f\in C(\bOm)$, since its proof via Dini's theorem uses only pointwise convergence. What is lost is the smoothing effect of Lemma~\ref{lem:modulus}, and hence local uniformity for non-continuous $f$. The same applies to averages against any radial probability density supported in $B_x$.
\end{remark}

\begin{question}
Does a geometric rate as in Theorem~\ref{thm:geometric} hold on every bounded $C^2$ (or $C^{1,1}$) domain when $\delta\ge\kappa\rho$? Near the boundary, $\rho^\gamma$ is strictly superharmonic for $\gamma\in(0,1)$, which suggests an affirmative answer; the difficulty is to construct a single Lyapunov function $V$ on $\Om$ with $\sigma V\le\lambda V$.
\end{question}

\begin{question}
What is the exact asymptotic rate of $\norm{\sigma^kw}$, for instance for $\delta=\rho$ on a ball? It is natural to expect $\norm{\sigma^kw}^{1/k}\to\lambda_*$ where $\lambda_*$ is a spectral quantity of $\sigma$ acting on a suitable weighted space.
\end{question}

\begin{question}
Theorem~HN also covers measurable $\delta$ that are locally bounded away from zero. To what extent do Theorems~\ref{thm:A}--\ref{thm:B} extend to discontinuous $\delta$, where $\sigma$ need not preserve continuity?
\end{question}


\begin{thebibliography}{99}

\bibitem{AG01}
D.~H. Armitage and S.~J. Gardiner, \emph{Classical Potential Theory}, Springer Monographs in Mathematics, Springer-Verlag, London, 2001.

\bibitem{Doo84}
J.~L. Doob, \emph{Classical Potential Theory and Its Probabilistic Counterpart}, Grundlehren der mathematischen Wissenschaften 262, Springer-Verlag, New York, 1984.

\bibitem{Eck93}
J.~Eckhoff, Helly, Radon, and Carath\'eodory type theorems, in: \emph{Handbook of Convex Geometry} (P.~M. Gruber and J.~M. Wills, eds.), Vol.~A, North-Holland, Amsterdam, 1993, pp.~389--448.

\bibitem{GT83}
D.~Gilbarg and N.~S. Trudinger, \emph{Elliptic Partial Differential Equations of Second Order}, 2nd ed., Springer-Verlag, Berlin, 1983.

\bibitem{Han96}
W.~Hansen, Restricted mean value property and harmonic functions, in: \emph{Potential Theory---ICPT 94} (Kouty, 1994), de Gruyter, Berlin, 1996, pp.~67--90.

\bibitem{HN93}
W.~Hansen and N.~Nadirashvili, A converse to the mean value theorem for harmonic functions, \emph{Acta Math.} \textbf{171} (1993), 139--163.

\bibitem{HN94}
W.~Hansen and N.~Nadirashvili, Littlewood's one circle problem, \emph{J. London Math. Soc. (2)} \textbf{50} (1994), 349--360.

\bibitem{Jav07}
M.~Javaheri, Harmonic functions via restricted mean-value theorems, preprint, \texttt{arXiv:0709.3311}, 2007.

\bibitem{Kel34}
O.~D. Kellogg, Converses of Gauss' theorem on the arithmetic mean, \emph{Trans. Amer. Math. Soc.} \textbf{36} (1934), 227--242.

\bibitem{Mul56}
M.~E. Muller, Some continuous Monte Carlo methods for the Dirichlet problem, \emph{Ann. Math. Statist.} \textbf{27} (1956), 569--589.

\bibitem{NV94}
I.~Netuka and J.~Vesel\'y, Mean value property and harmonic functions, in: \emph{Classical and Modern Potential Theory and Applications} (K.~GowriSankaran et al., eds.), NATO ASI Series C 430, Kluwer, Dordrecht, 1994, pp.~359--398.

\bibitem{Vee72}
W.~A. Veech, A converse to Gauss' theorem, \emph{Bull. Amer. Math. Soc.} \textbf{78} (1972), 444--446.

\bibitem{Vee73}
W.~A. Veech, A zero-one law for a class of random walks and a converse to Gauss' mean value theorem, \emph{Ann. of Math. (2)} \textbf{97} (1973), 189--216.

\bibitem{Vol09}
V.~Volterra, Alcune osservazioni sopra propriet\`a atte individuare una funzione, \emph{Atti Accad. Naz. Lincei Rend. (5)} \textbf{18} (1909), 263--266.

\end{thebibliography}
\end{document}